\documentclass[a4paper,12pt]{article}
\usepackage{comment}
\usepackage{cite}
\usepackage{amsmath}
\usepackage{amssymb}
\usepackage{amsfonts}
\usepackage[T1]{fontenc}
\usepackage[utf8]{inputenc}
\usepackage{graphicx}
\usepackage{fancyhdr}
\usepackage{float}
\usepackage{xcolor}
\usepackage{authblk}
\usepackage{mathrsfs}
\usepackage{empheq}

\usepackage{geometry}
\begin{document}

\title{Successive approximations of $\pi$ using Euler Beta functions}

\author[$\dagger$]{Jean-Christophe {\sc Pain}\\
\small
CEA, DAM, DIF, F-91297 Arpajon, France\\
Université Paris-Saclay, CEA, Laboratoire Matière en Conditions Extr\^emes,\\ 
91680 Bruyères-le-Ch\^atel, France
}

\maketitle

\begin{abstract}
In this didactic note, we describe a procedure to derive successive approximations of $\pi$ using Euler Beta functions. It is an interesting exercise for undergraduate students, since it involves polynomial roots, integral calculations, inequalities and Euler Beta functions.
\end{abstract}

\section{General procedure}

Let us consider the polynomials
\begin{equation}
\mathscr{P}_q(x)=4+\lambda_qx^{4q}(1-x)^{4q},
\end{equation}
where $q$ is a positive integer. The polynomial $x^2+1$ divides $\mathscr{P}_q(x)$ if and only if
\begin{equation}
\mathscr{P}_q(i)=\mathscr{P}_q(-i)=4+(-4)^q\lambda_q=0.
\end{equation}
This leads to the solution
\begin{equation}
\lambda_q=\left(-\frac{1}{4}\right)^{q-1}
\end{equation}
i.e. $\lambda_1=1$, $\lambda_2=-\frac{1}{4}$ and $\lambda_3=\frac{1}{16}$. Since
\begin{equation}
\int_0^1\frac{4}{x^2+1}dx=\pi,
\end{equation}
and defining
\begin{equation}
\mathscr{A}_q=\int_0^1\frac{\mathscr{P}_q(x)}{x^2+1}dx,
\end{equation}
we get
\begin{equation}
\mathscr{A}_q=\pi+\lambda_q\int_0^1\frac{x^{4q}(1-x)^{4q}}{x^2+1}dx.
\end{equation}
The fact that $x\in\left[0,1\right]$ enables one to write
\begin{equation}
\frac{x^{4q}(1-x)^{4q}}{2}<\frac{x^{4q}(1-x)^{4q}}{x^2+1}<x^{4q}(1-x)^{4q}.
\end{equation}
In the case where $q$ is even ($q=2p$), one has $\lambda_{2p}<0$, yielding
\begin{equation}
\frac{\lambda_{2p}}{2}\mathscr{K}_{2p}>\mathscr{A}_{2p}-\pi>\lambda_{2p}\mathscr{K}_{2p}
\end{equation}
or
\begin{equation}
\mathscr{A}_{2p}-\frac{\lambda_{2p}}{2}\mathscr{K}_{2p}<\pi<\mathscr{A}_{2p}-\lambda_{2p}\mathscr{K}_{2p},
\end{equation}
where
\begin{equation}
\mathscr{K}_{q}=\int_0^1x^{4q}(1-x)^{4q}dx=B(4q+1,4q+1)
\end{equation}
where 
\begin{equation}
B(r,s)=\int_0^1x^{r-1}(1-x)^{s-1}dx
\end{equation}
is the usual Euler Beta function \cite{Abra}. The integrand is continuous over $\left]0,1\right[$ and equivalent to $x^{r-1}$ in 0 and to $(1-x)^{s-1}$ in 1, therefore the integral converges if and only if $r,s>0$. One has
\begin{equation}
\Gamma(r)\,\Gamma(s)=\int_0^{+\infty}z^{r-1}e^{-z}\mathrm dz\int_0^{+\infty}t^{s-1}e^{-t}dt.
\end{equation}
Let us consider the change of variables $z=uv,~t=u\left(1-v\right)$, which Jacobian is 

\begin{equation}
\left|\begin{array}{cc}
\cfrac{\partial z}{\partial v} & \cfrac{\partial z}{\partial u}\\
\cfrac{\partial t}{\partial v} & \cfrac{\partial t}{\partial u}
\end{array}\right|=u.
\end{equation}
Then, the Fubini theorem gives
\begin{equation}
\Gamma(r)\,\Gamma(s)=\int_0^{+\infty}u^{r+s-1}\operatorname e^{-u}\,\mathrm du\int_0^1 v^{r-1}(1-v)^{s-1}\,\mathrm dv=\Gamma(r+s)\,\mathrm{B}(r,s). 
\end{equation}
If $r$ and $s$ are integers larger than 1, one has
\begin{equation}
B(r,s)=\frac{(r-1)!(s-1)!}{(r+s-1)!}.
\end{equation}
This enables us to write
\begin{empheq}[box=\fbox]{align}\label{even}
\mathscr{A}_{2p}+\frac{1}{2^{4p-1}}B(8p+1,8p+1)<\pi<\mathscr{A}_{2p}+\frac{1}{2^{4p-2}}B(8p+1,8p+1).
\end{empheq}
In the case where $q$ is odd ($q=2p+1$), one has $\lambda_{2p+1}>0$ which gives
\begin{equation}
\frac{\lambda_{2p+1}}{2}\mathscr{K}_{2p+1}<\mathscr{A}_{2p+1}-\pi<\lambda_{2p+1}\mathscr{K}_{2p+1},
\end{equation}
i.e.
\begin{equation}
\mathscr{A}_{2p+1}-\lambda_{2p+1}\mathscr{K}_{2p+1}<\pi<\mathscr{A}_{2p+1}-\frac{\lambda_{2p+1}}{2}\mathscr{K}_{2p+1}
\end{equation}
or
\begin{empheq}[box=\fbox]{align}\label{odd}
\mathscr{A}_{2p+1}-\frac{1}{2^{4p}}B(8p+5,8p+5)<\pi<\mathscr{A}_{2p+1}-\frac{1}{2^{4p+1}}B(8p+5,8p+5).
\end{empheq}

\section{Numerical values}

\begin{table}[ht!]
    \centering
    \begin{tabular}{c|c}\hline\hline
 $p$ & $\mathscr{A}_p$ \\\hline\hline
  1  & $\cfrac{22}{7}$ \\
     & \\
  2  & $\cfrac{47171}{15015}$ \\
     & \\
  3  & $\cfrac{431302721}{137287920}$ \\
     & \\
  4  & $\cfrac{741269838109}{235953517800}$ \\
     & \\
  5  & $\cfrac{26856502742629699}{8548690331301120}$ \\
     & \\ \hline\hline
    \end{tabular}
    \caption{Coefficients $\mathscr{A}_p$ for $p$=1,2,3, 4 and 5.}
    \label{tab1}
\end{table}

\begin{table}[ht!]
    \centering
    \begin{tabular}{c|c}\hline\hline
 $p$ & $B(8p+1,8p+1)$ \\\hline\hline
  1  & $\cfrac{1}{218790}$ \\
     & \\
  2  & $\cfrac{1}{19835652870}$ \\
     & \\
  3  & $\cfrac{1}{1580132580471900}$ \\
     & \\
  4  & $\cfrac{1}{119120569161268384710}$ \\
     & \\
  5  & $\cfrac{1}{8708083907400230293391220}$ \\  
     & \\ \hline\hline
    \end{tabular}
    \caption{Values of the Beta function $B(8p+1,8p+1)$ for $p$=1,2,3, 4 and 5.}
    \label{tab:my_label}
\end{table}

\begin{table}[ht!]
    \centering
    \begin{tabular}{c|c}\hline\hline
 $p$ & $B(8p+5,8p+5)$ \\\hline\hline
  1  & $\cfrac{1}{67603900}$ \\
     & \\
  2  & $\cfrac{1}{5651707681620}$ \\
     & \\
  3  & $\cfrac{1}{435975364243345080}$ \\
     & \\
  4  & $\cfrac{1}{32303415440209084881892}$ \\
     & \\
  5  & $\cfrac{1}{2336116978969951755817600200}$ \\
     & \\\hline\hline
    \end{tabular}
    \caption{Values of the Beta function $B(8p+5,8p+5)$ for $p$=1,2,3, 4 and 5.}
    \label{tab2}
\end{table}

\begin{table}[ht!]
\centering
\begin{tabular}{c|c|c}\hline\hline
 $q$ & Lower bound & Upper bound \\\hline\hline
  2 & 3.1415917420539930298982661346371132 & 3.1415917425162444680549405276824849 \\
  3 & 3.1415926534037147847889535869075511 & 3.1415926538659662229456279799529227\\
  4 & 3.1415926535891641685926362710915484 & 3.1415926535891645141740268925113291\\
  5 & 3.1415926535897930996904718896604452 & 3.1415926535897934452718625110802259\\
  6 & 3.1415926535897932379689068638929499 & 3.1415926535897932379691868574946798\\
  7 & 3.1415926535897932384625310710253762 & 3.1415926535897932384628110646271061\\
  8 & 3.1415926535897932384626429738632550 & 3.1415926535897932384626429740994341\\
  9 & 3.1415926535897932384626433831848238 & 3.1415926535897932384626433834210030\\
  10 & 3.1415926535897932384626433832791528 & 3.1415926535897932384626433832791530\\
  11 & 3.1415926535897932384626433832795028 & 3.1415926535897932384626433832795030\\\hline\hline
\end{tabular}
\caption{Values of the lower and upper bounds given by Eqs. (\ref{even}) and (\ref{odd}). For $q=2p$, $p=1, \cdots, 5$, the bounds are given by Eq. (\ref{even}) and for $q=2p+1$, $p=1, \cdots, 5$, the bounds are given by Eq. (\ref{odd}). For instance, in the case $q=2$, the lower bound is $\mathscr{A}_{2}+\frac{1}{2^{3}}B(9,9)$ (second column) and the upper bound $\mathscr{A}_{2}+\frac{1}{2^{2}}B(9,9)$ (third column). In the case $q=3$, the lower bound is $\mathscr{A}_{3}-\frac{1}{2^{4}}B(13,13)$ (second column) and the upper bound $\mathscr{A}_{3}-\frac{1}{2^{5}}B(13,13)$ (third column). We took 35 figures after the comma.}\label{tab3}
\end{table}

\clearpage

The exact value of $\pi$ with 40 significant digits after the comma is
\begin{equation}
\pi=3.141592653589793238462643383279502884197.
\end{equation}
One possibility to remember the first 126 first digits of $\pi$ is, in french, to remember the following text (attributed to Maurice Decerf) and count the number of letters of each word (when a word has 10 letters, it means that the digit is zero):

\vspace{5mm}

Que j’aime à faire apprendre un nombre utile aux sages !

Glorieux Archimède, artiste ingénieux,

Toi de qui Syracuse aime encore la gloire,

Soit ton nom conservé par de savants grimoires !

Jadis, mystérieux, un problème bloquait

Tout l’admirable procédé, l’œuvre grandiose

Que Pythagore découvrit aux anciens Grecs.

O, quadrature ! Vieux tourment du philosophe !

Insoluble rondeur, trop longtemps vous avez

Défié Pythagore et ses imitateurs.

Comment intégrer l’espace bien circulaire ?

Former un triangle auquel il équivaudra ?

Nouvelle invention : Archimède inscrira

Dedans un hexagone, appréciera son aire,

Fonction du rayon. Pas trop ne s’y tiendra

Dédoublera chaque élément antérieur ;

Toujours de l’orbe calculée approchera ;

Définira limite ; enfin, l’arc, le limiteur

De cet inquiétant cercle, ennemi trop rebelle !

Professeur, enseignez son problème avec zèle ! 

\section{Expression of integral $\mathscr{A}_q$}

The general value of integral $\mathscr{A}_q$ is \cite{Mathematica}:
\begin{eqnarray}
\mathscr{A}_q&=&\int_0^1\frac{\mathscr{P}_q(x)}{x^2+1}dx\nonumber\\
&=&\pi+\frac{(-1)^{q+1}}{2^{10q-1}}\sqrt{\pi}\frac{\Gamma(4q+1)}{\Gamma(4q+\frac{3}{2}}~_3F_2\left[
\begin{array}{c}
1, 2q+\frac{1}{2}, 2q+1\\
4q+1, 4q+\frac{3}{2}
\end{array};-1\right].
\end{eqnarray}
However, when $q$ is specified, the expression of $\mathscr{P}_q$ is simple, and the integral $\mathscr{A}_q$ can be calculated analytically. For instance, one has
\begin{equation}
\mathscr{P}_1(x)=4-4x^2+5x^4-4x^5+x^6,
\end{equation}
which integral between 0 and 1 is equal to $22/7$, as mentioned in Table 1. In the same way, one has
\begin{equation}
\mathscr{P}_2(x)=4-4x^2+4x^4-4x^6+\frac{15}{4}x^8+2x^ 9-\frac{43}{4}x^{10}+12x^ {11}-\frac{27}{4}x^ {12}+2x^ {13}-\frac{x^{14}}{4},
\end{equation}
which integral between 0 and 1 is equal to $47171/15015$, as mentioned in Table 1, and

\begin{eqnarray}
\mathscr{P}_3(x)&=&4-4x^2+4x^4-4x^6+4x^8-4x^{10}+\frac{65}{16}x^{12}-\frac{3}{4}x^{13}+\frac{x^{14}}{16}-13x^ {15}+\frac{247}{8}x^{16}-\frac{73}{2}x^{17}\nonumber\\
& &+\frac{215}{8}x^{18}-13x^{19}+\frac{65}{16}x^{20}-\frac{3}{4}x^{21}+\frac{x^{22}}{16},
\end{eqnarray}
which integral between 0 and 1 is equal to $431302721/137287920$.

\end{document}